\documentclass[12pt]{amsart}
\usepackage{amsthm}
\usepackage{amsmath}
\usepackage{amsfonts}
\newtheorem{theorem}{Theorem}[section]
\newtheorem{lemma}[theorem]{Lemma}
\newtheorem{proposition}[theorem]{Proposition}
\newtheorem{corollary}[theorem]{Corollary}
\newtheorem{conjecture}[theorem]{Conjecture}
\newtheorem{definition}[theorem]{Definition}
\newtheorem{example}[theorem]{Example}
\newtheorem{remark}[theorem]{Remark}

\def\KK{{\mathbb K}}

\def\AA{{\mathcal A}}

\def\OO{{\mathcal O}}

\def\MM{{\mathbb M}}
\def\ss{{\mathcal S}}
\def\aa{{\mathbb A}}

\def\mm{{\mathbb M}}
\def\CC{{\mathbb C}}
\def\PP{{\mathbb P}}
\def\NN{{\mathbb N}}
\def\QQ{{\mathbb Q}}
\def\RR{{\mathbb R}}
\def\ZZ{{\mathbb Z}}
\def\Res{\mbox{\bf Residue}\,}
\def\Del{{\Delta_{\bar{\sigma}}}}

\begin{document}
\title{Macaulay style formulas for toric residues}
\author{Carlos D' Andrea}
\address{Miller Institute for Basic Research in Science and Department of Mathematics
University of California at Berkeley, Berkeley CA 94720, USA}
\email{cdandrea@math.berkeley.edu}
\author{Amit Khetan}
\address{Department of Mathematics
University of California at Berkeley, Berkeley CA 94720, USA}
\email{akhetan@math.berkeley.edu}

\begin{abstract}
We present an explicit formula for computing toric residues as a
quotient of two determinants, \'a la Macaulay, where the numerator is a
minor of the denominator.  We also give an irreducible representation
of toric residues by extending the theory of subresultants to
monomials of critical degree in the homogeneous coordinate ring of the
corresponding toric variety.
\end{abstract}

\date{\today}
\maketitle

\section{introduction}

The toric residue of $n+1$ divisors on an $n$-dimensional toric
variety was first introduced by Cox \cite{cox} in the case when all
divisors are ample of the same class, and extended to the general case
by Cattani, Cox, and Dickenstein \cite{CCD}.  Toric residues have been
found to be useful in a variety of contexts such as mirror symmetry
\cite{BM}, the Hodge structure of hypersurfaces \cite{BC}, and in the
study of sparse resultants \cite{CDS}.

Another related, perhaps more familiar, notion is the global residue
in the torus.  Given a system $f_1, \dots, f_n$ of $n$ Laurent
polynomials in $n$ variables with a finite set of common zeroes in the
torus $T = (C^{\ast})^n,$ and another Laurent polynomial $q$, the
global residue of $q$ with respect to $f_1, \dots, f_n$ is the sum of
the Grothendieck residues of the differential form
$$\phi_q = \frac{q}{f_1 \cdots f_n} \frac{{\rm d}t_1}{t_1} \wedge
\cdots \wedge \frac{{\rm d}t_n}{t_n},$$

\noindent at each zero of the $f_i$. This turns out to be a rational
function in the coefficients in the $f_i$ and has a wide variety of
applications in algebra and analysis.  The residue in the torus has
been studied by Khovanskii, Gelfond, and Soprounov \cite{Kho,GK,
Sop}. Cattani and Dickenstein \cite[Theorem 4]{CD} showed that the
global residue in the torus is equal to a particular toric residue in
the sense above of \cite{cox, CCD}. The precise relationship between
the toric residue and the global residue in the torus is discussed
in Section \ref{smac}.

In this paper, we present an explicit formula for computing a toric
residue as a quotient of two determinants, where the numerator is a
minor of the denominator. This is an improvement over earlier
algorithms in \cite{CCD,CD} by eliminating costly Gr\"obner basis
computations.

Both the numerator and denominator of our quotient formula turn out to
be divisible by the same extraneous factor. It would be useful to have a
description of the residue in reduced form. Indeed, the denominator of
the toric residue has already been identified with the sparse
resultant \cite{CDS}. Our second main result is an identification of
the numerator with a ``toric subresultant'', analagous to the
multivariate subresultant of Chardin \cite{Cha1}.
In the dense case, this numerator and its properties have been deeply studied by Jouanolou in \cite{J1,J2,J3}.
Our results may be regarded as a generalization of Jouanolou's work.

We start with some notation on toric varieties. For more details we
refer the reader to \cite{ful,cox1,cox}.  Let $X$ be a projective
toric variety of dimension $n$, hence determined by a fan $\Sigma
\subset \mathbb{R}^n$. The generators of the $1$-dimensional cones in
$\Sigma$ will be denoted $\eta_0, \dots, \eta_{s-1}$. The Chow group
$A_{n-1}(X)$ has rank $s-n$. We work in the polynomial ring $S :=
\aa[x_0, \dots, x_{s-1}]$ where each variable $x_i$ corresponds to ray
$\eta_i$ and hence to a torus-invariant divisor $D_i$ of $X$. We grade
$S$ by declaring that the monomial $\prod x_i^{a_i}$ has degree $[\sum
a_i D_i] \in A_{n-1}(X)$. The base ring $\aa$ will be the coefficient
space of our polynomial system and is specified below. We abbreviate
$\beta_0 := [\sum_iD_i]$, the anticanonical class on $X$.  The {\em
irrelevant ideal} $B(\Sigma)$ is generated by the elements
$\hat{x}_{\sigma} = \prod_{\eta_i \notin \sigma} x_i$ where $\sigma$
ranges over all $n$-dimensional cones in $\Sigma$.

We now recall the definition of the toric residue from \cite{cox}.
This depends on $n+1$ generic $S$-homogeneous polynomials which are
homogenizations of Laurent polynomials in $n$ variables.  Given a
polynomial of a certain critical degree we construct a differential
$n$-form which gives rise to a top cohomology class of the canonical
sheaf of differentials $\Omega_X^n = \OO(-\beta_0)$. The toric residue
is defined to be the trace of this cohomology class.

Formally, pick ample degrees $\alpha_0, \dots, \alpha_n$ and consider
generic polynomials:
\begin{equation}\label{generic}
F_i(u,x):=\sum_{a\in\AA_{i}}u_{ia}x^a,\ i=0,\dots,n.
\end{equation}
where $\AA_i:=\{a\in\NN^s:\,\deg(a)=\alpha_i \}.$ Set
$\aa:=\QQ[u_{ia};i=0,\dots,n;a\in\AA_i]$, and write $Q(\aa)$ for the
field of quotients of $\aa$.  Let $\rho = \sum_i\alpha_i -\beta_0$,
which is the {\em critical degree} of our system.  For any subset $I =
\{i_1, \dots, i_n\}$ of $\{0, \dots, s-1\}$ we write
$$ \det(\eta_I) := \det(\langle e_l, \eta_{i_j} \rangle_{1 \leq l, j
\leq n}), \quad dx_I = dx_{i_1} \wedge \cdots \wedge dx_{i_n}, \quad
\hat{x_I} = \prod_{j \notin I} x_j. $$
The {\em Euler form} on $X$ is the sum over all $n$-element subsets $I$
of $\{1, \dots, s\}$:
$$ \Omega := \sum_{|I| = n} \det(\eta_I) \hat{x_I} dx_I. $$ The
polynomials $F_i$ determine an open cover $U_i = \{x \in X \ : \
F_i(x) \neq 0 \}$. A polynomial $g \in S_{\rho}$ gives an element
$\omega_g = \frac{g \Omega}{F_0 \cdots F_n}$ which is a \u{C}ech
cocycle of degree $n$ with respect to the open cover $U_i$. Therefore,
there is an induced cohomology class $[\omega_g] \in H^n(X,
\Omega_X^n)$, and we define $\Res_F\,(g) = {\rm Tr}_X([\omega_g])$.
The toric residue is thus a map $\Res_F : S_{\rho} \to Q(\aa)$ .

The remainder of this paper is organized as follows: Section
\ref{delta0} gives a combinatorial construction of explicit elements
of $S_{\rho}$ with residue $\pm 1$, playing the role of the toric
Jacobian from \cite{CDS} and generalizing the elements
$\Delta_{\sigma}$ constructed in \cite{CCD} on simplicial toric
varieties.  This is used in Section \ref{mac} to give the promised
Macaulay style residue formula.

Sections 4 and 5 discuss the reduced numerator of the toric residue of
a monomial $h$. This is shown to be a toric subresultant associated to
the systen $F_0, \dots, F_n$ and the monomial $h$.  In the dense case,
the subresultant was introduced by \cite{Cha1}, and different
algorithms for computing subresultants have been developed in
\cite{Cha2,sza}.

In Section \ref{smac}, we will show how our results can be used to
compute global residues in the torus. We see how our results
generalize some explicit formulas given by Macaulay in \cite{Mac} for
computing global residues of dense homogeneous systems.

\section{Elements with nonzero residue} \label{delta0}

Let $X$ be a projective, but not necessarily simplicial, toric variety
with fan $\Sigma$ and $F_0, \dots, F_n$ generic elements of ample
degrees $\alpha_i$ as above.

 Pick a complete flag $\bar{\sigma} = \{0\} \subset \sigma_1 \subset
\cdots \subset \sigma_n$ where each $\sigma_i$ is a cone of dimension
$i$ in $\Sigma$.  For $i = 1, \dots, n$, let $z_i$ be the product of
all variables $x_j$ such that $\eta_j \in \sigma_i$ but $\eta_j \notin
\sigma_{i-1}$. We set $z_{n+1} = \prod_{\eta_j \notin \sigma_n} x_j$.

We will see that each $F_j$ can be written in the form
$$F_j = \sum_{i=1}^{n+1} A_{ij} z_i. $$
The $(n+1) \times (n+1)$-determinant $\Del = \det(A_{ij})$ is in
$S_{\rho}$. We have the following theorem generalizing
\cite[Theorem 0.2]{CCD}.

\begin{theorem} Suppose $X$, $F$, and $\bar{\sigma}$ are as above. Then
$\Res_F\,(\Del) \\ = \pm 1.$
\end{theorem}

Note that if $X$ were simplicial then each of $z_1, \dots, z_n$ would
be a single variable corresponding to the generators of the cone
$\sigma_n$ and $z_{n+1}$ would be the product of all of the remaining
variables. The element $\Del$, in this case, is the same as the
element $\Delta_{\sigma_n}$ from \cite{CCD}.

\begin{proof}

The $\eta_i$ are ample classes, therefore $\langle F_0, \dots, F_n
\rangle \subset B(\Sigma)$.  We next show that $B(\Sigma) \subset
\langle z_1, \dots z_{n+1} \rangle$.  Consider a generator
$\hat{x}_{\tau}$ where $\tau$ is a maximal cone of $\Sigma$. Recall
that $\hat{x}_{\tau} = \prod_{\eta_i \notin \tau} x_i$. Now $\tau \cap
\sigma_n$ is a face of $\sigma_n$. Choose $i$ such that $\tau \cap
\sigma_n \supset \sigma_{i-1}$ but $\tau \cap \sigma_n \not \subset
\sigma_{i}$. Now, we see that the one dimensional cones in $\tau \cap
\sigma_i$ are all contained in $\sigma_{i-1}$, and so none of the one
dimensional cones in $z_i$ are in $\tau$.  Therefore $z_i$ divides
$\hat{x}_{\tau}$ as desired. Hence we can write (nonuniquely) for $j =
0, \dots n$
$$F_j = \sum_{i=1}^{n+1} A_{ij} z_i. $$
Now an application of the global transformation law (Theorem 0.1 in
\cite{CCD}), shows that $\Res_F\,(\Del) = \Res_z\,(1)$. And so we need only prove that $\Res_z\,(1) = \pm 1$.

Let $X'$ be a new toric variety arising from a simplicial refinement
$\Sigma'$ of $\Sigma$ with the same $1$-dimensional cones, hence the
same coordinate ring, albeit with a smaller irrelevant ideal
$B(\Sigma') \subset B(\Sigma)$. We have a natural map $f : X' \to X$.
Let $\mathcal{U}$ be the open cover on $X$ defined by $\{ U_i = \{x
\in X \ : \ z_i \neq 0 \} \}$. We have an analagous collection of open
sets $\mathcal{U'}$ on $X'$ defined by the same equations $\{ U'_i = \{x \in X'
\ : \ z_i \neq 0 \} \}$. As $B(\Sigma') \subset B(\Sigma) \subset
\langle z_1, \dots z_{n+1} \rangle$, $\mathcal{U'}$ is also an open cover of
$X'$.

Now, on $X$, the element $1 \in S$ gives the $n$-form

\begin{equation} \label{omega}
\omega = \frac{\Omega}{z_1 \cdots z_{n+1}} = \frac{\Omega}{x_0 \cdots
x_{s-1}}.
\end{equation}

\noindent This form is defined on the open set $$U = \{x \in X \ : \
z_1(x)\neq0, \dots, z_{n+1}(x) \neq 0 \},$$ which is the same as the
open set defined by $x_0\neq0, \dots, x_{s-1} \neq 0$. The map $f$
defines an isomorphism on this open set, and we can pull back $\omega$
to get a new form $\omega' = f^{\ast}(\omega)$ on $f^{-1}(U) \subset
X'$, defined by the same formula (\ref{omega}).  Now,
$f^{-1}\mathcal{U}$ is an open cover of $X'$, and moreover it refines
the open cover $\mathcal{U}'$, since $z_i \neq 0$ on each
$f^{-1}(U_i)$. Since the natural map from \u{C}ech cohomology to sheaf
cohomology respects refinement $[\omega'_{\mathcal{U}'}] =
[\omega'_{f^{-1}\mathcal{U}}]$ as cohomology classes in $H^n(X',
\Omega_X'^n)$. Moreover, $f$ is birational so $Tr_X' \circ f^{\ast} =
Tr_X$.  Putting it together we have:

\begin{equation} \label{reseq}
 \Res_z(1)_X = \pm 1 \iff \Res_z(1)_{X'} = \pm 1.
\end{equation}

In \cite{CCD} it was shown that $Tr([\omega'_{\mathcal{V}}]) = \pm 1$
where $\mathcal{V}$ is an open cover coming from the simplicial
construction of $\Delta_{\sigma'}$, and $\omega'$ is the same $n$-form
as in (\ref{omega}). So, we need only show that the cohomology classes
$[\omega'_{\mathcal{U}'}]$ and $[\omega'_{\mathcal{V}}]$ coincide for
some top dimensional cone $\sigma'$ giving rise to the cover
$\mathcal{V}$.  This will be done by picking the appropriate
$\sigma'$, and thereby $\mathcal{V}$, and finding an open cover
$\mathcal{W}$ refining both $\mathcal{U}'$ and $\mathcal{V}$ for which
$\omega'$ still determines a \u{C}ech cocycle.

\begin{lemma}
There is a unique (simplicial) cone $\sigma' \in \Sigma'$ of dimension
$n$ generated by $\{ \eta_{i_1}, \dots \eta_{i_n} \}$ such that the
corresponding variables $x_{i_1}, \dots , x_{i_n}$ satisy $x_{i_k} \mid z_k$ and
moreover $z_1 = \dots =z_n = 0$ on $X'$ if and only if $x_{i_1} = \dots =
x_{i_n} = 0$.
\end{lemma}

\begin{proof} We proceed by induction to show that $z_1 = \dots= z_k = 0$ if and only if $x_{i_1} = \dots= x_{i_k} = 0$ 
for a unique $k$ dimensional cone $\sigma'_k \in \Sigma'$.  The base
case $k=1$ is trivial as $z_1$ is a single variable $x_{i_1}$
corresponding to a unique $1$ dimensional cone. For the inductive step
suppose $x_1 = \dots = x_{i_k} = 0$ and $z_{k+1} = 0$. Let $\sigma'_k$
be the cone corresponding to $x_{i_1}, \dots, x_{i_k}$. The torus
orbits of $X'$ are in one to one correspondance with the cones in
$\Sigma'$.  The orbit corresponding to a cone $\sigma'$ is specified
exactly by the vanishing of the variables corresponding to the
1-dimensional generators. In particular if some variable $x_{i_{k+1}}$
dividing $z_{k+1}$ vanishes, then there is a cone in $\Sigma'$
containing the rays corresponding to $x_{i_1}, \dots, x_{i_{k+1}}$. As
$\Sigma'$ is simplicial, the above set of rays must itself be a
$(k+1)$-dimensional cone in $\Sigma'$.  So, it is enough to show that
there is a unique such cone.

Recall that we have a complete flag of cones $\sigma_0 \subset \cdots
\subset \sigma_k \subset \cdots \sigma_n \subset \Sigma$. By
construction, $\sigma'_k \subset \sigma_k$ and $z_{k+1}$, is made up
of the rays in $\sigma_{k+1}$ but not in $\sigma_k$. Now $\sigma_k$
was a facet of $\sigma_{k+1}$ and $\sigma'_k$ is a cone in its
triangulation in $\Sigma'$.  Therefore, there is a unique cone
$\sigma'_{k+1}$ in the triangulation of $\sigma_{k+1}$ containing
$\sigma'_k$ (No two cones in a triangulation can meet in a facet of
the original cone). Let $x_{i_{k+1}}$ be the additional generator of
$\sigma'_{k+1}$. This completes the induction and we take $\sigma' =
\sigma'_n$.

\end{proof}

Now, let $z'_k = x_{i_k}$ for $k = 1, \dots, n$, and $z'_{n+1} =
\prod_{\eta_j \notin \sigma'} x_j$ the product of the remaining
variables. We define the open cover $\mathcal{V}$ by $$V_i = \{ x \in
X' \ : \ z'_i \neq 0 \}.$$ By the above lemma we have $U'_i \subset
V_i$ for $i = 1, \dots, n$, $V'_{i+1} \subset U'_{i+1}$, and
$\bigcup_{i=1}^n V_i = \bigcup_{i=1}^n U'_i$. Therefore the open cover
$\mathcal{W} = \{ U'_1, \dots, U'_n, V_{n+1} \}$ refines both
$\mathcal{U}'$ and $\mathcal{V}$, and both $\omega'_{\mathcal{U}'}$
and $\omega'_{\mathcal{V}}$ map to a well defined \u{C}ech cocycle
$\omega'_{\mathcal{W}}$.  Hence $\omega'_{\mathcal{U}'}$ and
$\omega'_{\mathcal{V}}$ induce the same cohomology class, and hence
have the same trace, as desired. So $\Res_x(1)_{X'} = \pm 1$ which completes the proof by virtue of (\ref{reseq}).  
\end{proof}

\begin{example}\label{octh}

Let $O \subset \RR^3$ be the octahedron with vertices $(\pm 1, 0 ,
0)$, $(0, \pm 1, 0)$, $(0, 0, \pm 1)$. This determines a toric variety
$X_O$ whose normal fan has generators \ $\eta_0 = (-1, -1, -1),$
$\eta_1 = (-1, -1, 1),$ $\eta_2 = (-1, 1, -1),$ $\dots,$ $\eta_7 = (1,
1, 1)$ and top dimensional cones determined by the spans of the sets
of four rays with one coordinate fixed. Pick the ample degrees
$\alpha_0 = \cdots = \alpha_3 = \sum_{i=0}^7 D_i$ and consider the
generic system:
$$
\begin{array}{lcl}
F_0&:=&a_0x_0^2x_1^2x_2^2x_3^2 + a_1 x_0^2x_1^2x_4^2x_5^2 + a_2 x_0^2x_2^2x_4^2x_6^2 + a_3x_0x_1x_2x_3x_4x_5x_6x_7 \\
&+& a_4 x_1^2x_3^2x_5^2x_7^2 + a_5 x_2^2x_3^2x_6^2x_7^2 + a_6x_4^2x_5^2x_6^2x_7^2 \\
F_1&:=&b_0x_0^2x_1^2x_2^2x_3^2 + b_1 x_0^2x_1^2x_4^2x_5^2 + b_2 x_0^2x_2^2x_4^2x_6^2 + b_3x_0x_1x_2x_3x_4x_5x_6x_7 \\
&+& b_4 x_1^2x_3^2x_5^2x_7^2 + b_5 x_2^2x_3^2x_6^2x_7^2 + b_6x_4^2x_5^2x_6^2x_7^2 \\
F_2&:=&c_0x_0^2x_1^2x_2^2x_3^2 + c_1 x_0^2x_1^2x_4^2x_5^2 + c_2 x_0^2x_2^2x_4^2x_6^2 + c_3x_0x_1x_2x_3x_4x_5x_6x_7 \\
&+& c_4 x_1^2x_3^2x_5^2x_7^2 + c_5 x_2^2x_3^2x_6^2x_7^2 + c_6x_4^2x_5^2x_6^2x_7^2 \\
F_3&:=&d_0x_0^2x_1^2x_2^2x_3^2 + d_1 x_0^2x_1^2x_4^2x_5^2 + d_2 x_0^2x_2^2x_4^2x_6^2 + d_3x_0x_1x_2x_3x_4x_5x_6x_7 \\
&+& d_4 x_1^2x_3^2x_5^2x_7^2 + d_5 x_2^2x_3^2x_6^2x_7^2 + d_6x_4^2x_5^2x_6^2x_7^2
\end{array}$$

\noindent As our complete flag $\bar{\sigma}$ we pick $\{0\} \subset \{\eta_0\} \subset \{\eta_0, \eta_1\} \subset \{\eta_0, \eta_1, \eta_2, \eta_3\}$. This gives $z_0 = x_0,\  z_1 = x_1, \ z_2 = x_2x_3, \ z_3 = x_4x_5x_6x_7$ and we can write
$$
\begin{array}{lcl}
F_0 &=& z_0(a_0x_0x_1^2x_2^2x_3^2 + a_1x_0x_1^2x_4^2x_5^2 + a_2x_0x_2^2x_4^2x_6^2 + a_3x_1x_2x_3x_4x_5x_6x_7)\\
&+& z_1(a_4x_1x_3^2x_5^2x_7^2) + z_2(a_5x_2x_3x_6^2x_7^2) + z_3(a_6x_4x_5x_6x_7),
\end{array}$$

\noindent and similarly for $F_1, F_2, F_3$. Therefore:
\begin{eqnarray*}
\Del &=& \det \left(
{\begin{array}{cccc}
a_0x_0x_1^2x_2^2x_3^2 + \cdots & a_4x_1x_3^2x_5^2x_7^2 &
a_5x_2x_3x_6^2x_7^2 & a_6x_4x_5x_6x_7 \\
b_0x_0x_1^2x_2^2x_3^2 + \cdots & b_4x_1x_3^2x_5^2x_7^2 &
b_5x_2x_3x_6^2x_7^2 & b_6x_4x_5x_6x_7 \\
c_0x_0x_1^2x_2^2x_3^2 + \cdots & c_4x_1x_3^2x_5^2x_7^2 &
c_5x_2x_3x_6^2x_7^2 & c_6x_4x_5x_6x_7 \\
d_0x_0x_1^2x_2^2x_3^2 + \cdots & d_4x_1x_3^2x_5^2x_7^2 &
d_5x_2x_3x_6^2x_7^2 & d_6x_4x_5x_6x_7 \\
\end{array}}
\right) \\
&=& [0456]x_0x_1^3x_2^3x_3^5x_4x_5^3x_6^3x_7^5 + [1456]x_0x_1^3x_2x_3^3x_4^3x_5^5x_6^3x_7^5 \\
&+& [2456]x_0x_1x_2^3x_3^3x_4^3x_5^3x_6^5x_7^5 + [3456]x_1^2x_2^2x_3^4x_4^2x_5^4x_6^4x_7^6.
\end{eqnarray*}

\noindent Here the ``bracket'' $[0456]$ denotes the $4 \times 4$
determinant:
$$
\det \left(
{\begin{array}{cccc}
a_0 & a_4 & a_5 & a_6 \\
b_0 & b_4 & b_5 & b_6 \\
c_0 & c_4 & c_5 & c_6 \\
d_0 & d_4 & d_5 & d_6 \\
\end{array}}
\right).
$$

\end{example}

\begin{remark}
The construction uses that the $\alpha_i$ are ample degrees. It is
still an open problem to find an explicit element of nonzero residue
in a more general setting, for example when the $\alpha_i$ correspond to
nef and big divisors, i.e. have $n$-dimensional support. Indeed this is the
only obstruction to generalizing the Macaulay style formula of the next
section to this more general case.
\end{remark}

\section{Macaulay Style Formulas for Residues}\label{mac}

We now show how to use the element $\Del$ to give an explicit Macaulay
formula for the residue. We will need the following result.
Consider the map of free $\aa$ modules:

\begin{equation}\label{modc}
\begin{array}{cccc}
\phi:&S_{\rho-\alpha_0}\oplus \dots \oplus S_{\rho-\alpha_n}\oplus\aa
&\to&S_\rho\\ &(G_0,\dots,G_n,c)&\mapsto&\sum_{i=0}^n G_iF_i+c\,\Del.
\end{array}
\end{equation}

\begin{theorem} \label{codimthm} (Codimension 1 theorem)
Let $X$ be a projective toric variety, and $F$ a generic ample system
as above, then the map $\phi$ described above is generically
surjective. Equivalently, the degree $\rho$ component of the quotient
$S_F:=S / \langle F_0,\dots,F_n\rangle$ has $Q(\aa)$ dimension $1.$
\end{theorem}

The proof is postponed to the next section.  When $X$ is
simplicial, this result is due to \cite{CCD}. In a forthcoming
paper, Cox and Dickenstein \cite{CoD} prove a much more general codimension
theorem which implies the theorem above.

Let $\mm$ be the matrix associated with the $\aa$-linear map $\phi$ in
the monomial bases. As in \cite{Mac}, we shall index the rows of $\mm$
in correspondence with the elements of the monomial basis of the domain. Fix a monomial
$h\in S_\rho,$ and let $\tilde\mm$ be any square maximal submatrix
with nonvanishing determinant.  It turns out that one of the rows of
$\tilde\mm$ must be indexed by $(0,\dots,0,1).$ This is due to the
fact that $\Del$ does not belong to the ideal generated by
$F_0,\dots,F_n$ if this family does not have any common zero in
$X$. Let $\tilde \mm_h$ be the square submatrix of $\tilde\mm$ made by
deleting the row indexed by $(0,\dots,0,1)$ and the column indexed by
$h$. Now we are ready for the main result of this section.

\par
\begin{theorem}\label{quot}
$$\Res_F\,(h)=\pm
\frac{\det(\tilde\mm_h)}{\det(\tilde\mm)}.$$
\end{theorem}
\begin{proof}

Let $\tilde{M}^h$ be the matrix $\tilde\mm$  modified as follows: we multiply by $h$ all the elements in the column indexed by $h.$

\par Then, it turns out that $\det(\tilde{M}^h)=h\det(\tilde\mm).$
On the other hand, performing elementary operations in the columns of
$\tilde\mm$ and expanding the determinant along the column indexed by $h$,
it turns out that
$$h\det(\tilde\mm)=
\det(\tilde{M}^h)=G_0F_0+\dots+G_nF_n\pm\det(\tilde\mm_h)\,\Del.$$
By taking the residue of both sides, we get
$$\Res_F\,(h)\det(\tilde\mm)=\pm \det(\tilde\mm_h) \Res_F\,(\Del) = \pm \det (\tilde\mm_h)$$

\end{proof}
\begin{remark}
In the case of ``generalized unmixed systems'', when each $\alpha_i$
is an integer multiple of a fixed ample degree $\alpha$, one can also
use the toric Jacobian $J(F)$, scaled by an appropriate constant, in
place of the element $\Del$. The Jacobian has the computational
disadvantage of having much larger support than $\Del$ although it has
the advantage of being intrinsic to the toric system and not dependent
on a choice of $\bar{\sigma}$.

\end{remark}
\medskip
The following result is a straightforward consequence of Theorem \ref{quot},
and says that for \textit{any} polynomial $P$ of critical degree, $\Res_F\,(P)$
may be computed as a quotient of two determinants.
\begin{corollary}\label{gral}
Let $P=\sum_{\deg(a)=\rho}p_ax^a,$ where $p_a$ are constants. As before,
let  $\tilde\mm$ any square maximal submatrix
of $\mm$ having nonzero determinant.
Let $\tilde\mm_P$ be the matrix $\tilde\mm$ modified as follows:
for every monomial $a$ of critical degree, we replace the input indexed by
$((0,\dots,0,1);a)$  with the coefficient $p_a.$
Then, it turns out that
$$\Res_F\,(P)=  \pm
\frac{\det(\tilde\mm_P)}{\det(\tilde\mm)}.$$
\end{corollary}

\begin{proof}
Expand $\det(\tilde\mm_P)$ by the row indexed by
$\Del$, and use the linearity of the residue and
Theorem \ref{quot}.
\end{proof}

\begin{example}
Let  $X=\PP^1\times\PP^1$ whose fan has the following $1$-dimensional generators: $\eta_1=(1,0),\,\eta_2=(0,-1),\,\eta_3=(-1,0)$ and $\eta_4=(0,1).$
We pick the ample degrees $\alpha_0=D_2+D_3,\,\alpha_1=2D_2+D_3,\,\alpha_3=D_2+2D_3$ and consider the following generic polynomials having those degrees:
$$\begin{array}{lcl}
F_0&:=&a_0x_2x_3+a_1x_1x_2+a_2x_1x_4+a_3x_3x_4,\\
F_1&:=&b_0x_2^2x_3+b_1x_1x_2^2+b_2x_2x_2x_4+b_3x_1x_2x_4+b_4x_3x_4^2+b_5x_1x_4^2,\\
F_2&:=&c_0x_2x_3^2+c_1x_1x_2x_3+c_2x_1^2x_2+c_3x_3^2x_4+c_4x_1x_3x_4+c_5x_1^2x_4.
\end{array}$$
The critical degree is $-D_1+2D_2+2D_3-D_4$ and may be identified with
the set of nine integer points lying in the interior of a $3\times3$
square having integer vertices and edges parallel to the axes (see
\cite{CDS}).  To compute $\Del$ we can take $z_1 = x_1$, $z_2 = x_2$,
$z_3 = x_3x_4$.
$$ \Delta_z = \det \left(
\begin{array}{ccc}
a_1x_2+a_2x_4 & a_0x_3 & a_3 \\
b_1x_2^2+b_3x_2x_4+b_5x_4^2 & b_0x_2x_3+b_2x_3x_4 & b_4x_4 \\
c_1x_2x_3+c_2x_1x_2+c_4x_3x_4+c_5x_1x_4 &c_0x_3^2 & c_3x_3
\end{array}
\right). $$
In this case, $\MM$ is the following $9\times9$ matrix:
$$\MM=\bordermatrix{\cr
&a_3&a_2&0&a_0&a_1&0&0&0&0\cr
&0&a_3&a_2&0&a_0&a_1&0&0&0\cr
&0&0&0&a_3&a_2&0&a_0&a_1&0 \cr
&0&0&0&0&a_3&a_2&0&a_0&a_1\cr
&b_4&b_5&0&b_2&b_3&0&b_0&b_1&0\cr
&0&b_4&b_5&0&b_2&b_3&0&b_0&b_1\cr
&0&0&0&c_3&c_4&c_5&c_0&c_1&c_2\cr
&c_3&c_4&c_5&c_0&c_1&c_2&0&0&0\cr
&{\bf d_1}&{\bf d_2}&{\bf d_3}&{\bf d_4}&{\bf d_5}&{\bf d_6}&{\bf d_7}&{\bf d_8}&{\bf d_9}
},$$
where the ${\bf d_k}$ are the coefficients of $\Del$:
$$\begin{array}{lcl}
{\bf d_1} &:=& -a_2 b_4 c_0 + b_5 a_3 c_0 + a_2 c_3 b_2 + c_4 a_0 b_4 - c_4 a_3 b_2 - b_5 a_0 c_3  \\
{\bf d_2} &:=&   c_5 a_0 b_4 - c_5 a_3 b_2 \\
{\bf d_3} &:=&0 \\
{\bf d_4} &:=&   a_1 c_3 b_2 - a_1 b_4 c_0 + a_2 c_3 b_0 - b_3 a_0 c_3 + b_3 a_3 c_0 + c_1 a_0 b_4 - c_1 a_3 b_2 - c_4 a_3 b_0\\
{\bf d_5} &:=& c_2 a_0 b_4 - c_2 a_3 b_2 - c_5 a_3 b_0\\
{\bf d_6} &:=& 0 \\
{\bf d_7} &:=& a_1 c_3 b_0 - a_0 b_1 c_3 + a_3 b_1 c_0 - c_1 a_3 b_0 \\
{\bf d_8} &:=& 0 \\
{\bf d_9} &:=& 0. \\
\end{array}$$
It turns out that $\det(\MM)$ equals the sparse resultant of the $F_i$'s.
Let
$$\begin{array}{lcl}
P&:=&{\bf p_1}x_3^2x_4^2+{\bf p_2}x_1x_3x_4^2+{\bf p_3}x_1^2x_4^2+{\bf p_4}x_2x_3^2x_4+{\bf p_5}x_1x_2x_3x_4\\
&&+{\bf p_6}x_1^2x_2x_4+{\bf p_7}x_2^2x_3^2+{\bf p_8}x_1x_2^2x_3+{\bf p_9}x_1^2x_2^2
\end{array}$$
be any polynomial of critical degree. Following the notation of Corollary \ref{gral}, we have that $\tilde{\MM}=\MM$ and
$$\tilde{\MM_P}=\left(\begin{array}{ccccccccc}
a_3&a_2&0&a_0&a_1&0&0&0&0\\
0&a_3&a_2&0&a_0&a_1&0&0&0\\
0&0&0&a_3&a_2&0&a_0&a_1&0 \\
0&0&0&0&a_3&a_2&0&a_0&a_1\\
b_4&b_5&0&b_2&b_3&0&b_0&b_1&0\\
0&b_4&b_5&0&b_2&b_3&0&b_0&b_1\\
0&0&0&c_3&c_4&c_5&c_0&c_1&c_2\\
c_3&c_4&c_5&c_0&c_1&c_2&0&0&0\\
{\bf p_1}&{\bf p_2}&{\bf p_3}&{\bf p_4}&{\bf p_5}&{\bf p_6}&{\bf p_7}&{\bf p_8}&{\bf p_9}
\end{array}
\right).$$
So, we have that $\Res_F\,(P)=  \pm
\frac{\det(\tilde\mm_P)}{\det(\tilde\mm)}.$
\end{example}

\medskip
\begin{example}
Consider the generic system introduced in Example \ref{octh}. In this
case, $\MM$ is a $101\times63$ size matrix.  With the aid of Maple, we
have found a square maximal minor $\tilde\MM$ by choosing the rows
indexed by the element $\Del$, as constructed in Example $\ref{octh}$,
and the following monomials:
\begin{itemize}
\item In $S_{\rho-\alpha_0}:$
$$\begin{array}{l}
x_0^4 x_1^4 x_2^4 x_3^4,\, x_0^4 x_1^4 x_4^4 x_5^4,\, x_0^4 x_2^4 x_4^4 x_6^4,\,
x_1^4 x_3^4 x_5^4 x_7^4,\, x_2^4 x_3^4 x_6^4 x_7^4,\, x_4^4 x_5^4 x_6^4 x_7^4,\\
x_0^4 x_1^4 x_2^2 x_3^2 x_4^2 x_5^2,\, x_0^4 x_1^2 x_2^4 x_3^2 x_4^2 x_6^2,\,
x_0^3 x_1^3 x_2^3 x_3^3 x_4 x_5 x_6 x_7,\,
x_0^2 x_1^4 x_2^2 x_3^4 x_5^2 x_7^2,\\ x_0^2 x_1^2 x_2^4 x_3^4 x_6^2 x_7^2,\,
x_0^2 x_1^2 x_2^2 x_3^2 x_4^2 x_5^2 x_6^2 x_7^2,\,
x_0^4 x_1^2 x_4^4 x_5^2 x_2^2 x_6^2,\,
x_0^3 x_1^3 x_4^3 x_5^3 x_2 x_3 x_6 x_7,\\
x_0^2 x_1^4 x_4^2 x_5^4 x_3^2 x_7^2 ,\, x_0^2 x_1^2 x_4^4 x_5^4 x_6^2 x_7^2 ,\,
x_0^3 x_2^3 x_4^3 x_6^3  x_1 x_3 x_5 x_7,\,
x_0^2 x_2^4 x_4^2 x_6^4 x_3^2 x_7^2,\\ x_0^2 x_2^2 x_4^4 x_6^4 x_5^2 x_7^2,\,
x_0 x_1^3  x_2 x_3^3  x_4 x_5^3  x_6 x_7^3,\,
x_0 x_1 x_2^3  x_3^3  x_4 x_5 x_6^3  x_7^3,\\
x_0 x_1 x_2 x_3 x_4^3  x_5^3  x_6^3  x_7^3 ,\,
x_1^2  x_3^4  x_5^2  x_7^4  x_2^2  x_6^2 ,\, x_1^2  x_3^2  x_5^4  x_7^4  x_4^2  x_6^2 ,\,
x_2^2  x_3^2  x_6^4  x_7^4  x_4^2  x_5^2.
\end{array}$$
\item In $S_{\rho-\alpha_1}:$
$$\begin{array}{l}
x_0^4  x_1^4  x_2^4  x_3^4 ,\, x_0^4  x_1^4  x_4^4  x_5^4 ,\, x_0^4  x_2^4  x_4^4  x_6^4 ,\,
x_1^4  x_3^4  x_5^4  x_7^4 ,\, x_2^4  x_3^4  x_6^4  x_7^4 ,\, x_4^4  x_5^4  x_6^4  x_7^4 ,\\
x_0 ^4 x_1^4  x_2^2  x_3^2  x_4^2  x_5^2 ,\, x_0^4  x_1^2  x_2^4  x_3^2  x_4^2  x_6^2 ,\,
x_0^3  x_1^3  x_2^3  x_3^3  x_4 x_5 x_6 x_7,\,
x_0^2  x_1^4  x_2^2  x_3^4  x_5^2  x_7^2 ,\\ x_0^2  x_1^2  x_2^4  x_3^4  x_6^2  x_7^2 ,\,
x_0^4  x_1^2  x_4^4  x_5^2  x_2^2  x_6^2 ,\,
x_0^3  x_1^3  x_4^3  x_5^3  x_2 x_3 x_6 x_7,\,
x_0^2  x_1^4  x_4^2  x_5^4  x_3^2  x_7^2 ,\\
x_0^3 x_2^3  x_4^3  x_6^3  x_1 x_3 x_5 x_7,\,
x_0^2  x_2^4  x_4^2  x_6^4  x_3^2  x_7^2 ,\,
x_0 x_1^3  x_2 x_3^3  x_4 x_5^3  x_6 x_7^3 ,\,
x_0 x_1 x_2^3  x_3^3  x_4 x_5 x_6^3  x_7^3.
\end{array}$$

\item In $S_{\rho-\alpha_2}:$
$$\begin{array}{l}
x_0^4  x_1^4  x_2^4  x_3^4 ,\, x_0^4  x_1^4  x_4^4  x_5^4 ,\, x_0^4  x_2^4  x_4^4  x_6^4 ,\,
x_1^4  x_3^4  x_5^4  x_7^4 ,\, x_2^4  x_3^4  x_6^4  x_7^4 ,\, x_4^4  x_5^4  x_6^4  x_7^4 ,\\
x_0^4  x_1^4  x_2^2  x_3^2  x_4^2  x_5^2 ,\, x_0^4  x_1^2  x_2^4  x_3^2  x_4^2  x_6^2 ,\,
x_0^3  x_1^3  x_2^3  x_3^3  x_4 x_5 x_6 x_7,
x_0^2  x_1^4  x_2^2  x_3^4  x_5^2  x_7^2 ,\\ x_0^4  x_1^2  x_4^4  x_5^2  x_2^2  x_6^2 ,\,
x_0^3  x_1^3  x_4^3  x_5^3  x_2 x_3 x_6 x_7.
\end{array}$$
\item In $S_{\rho-\alpha_3}:$
$$\begin{array}{l}
x_0^4   x_1^4   x_2^4   x_3^4  ,\, x_0^4   x_1^4   x_4^4   x_5^4  ,\, x_0^4   x_2^2   x_4^4   x_6^4  ,\,
x_1^4   x_3^4   x_5^4   x_7^4  ,\, x_2^4   x_3^4   x_6^4   x_7^4  ,\\ x_4^4   x_5^4   x_6^4   x_7^4  ,\,
x_0^4   x_1^4   x_2^2   x_3^2   x_4^2   x_5^2.
\end{array}
$$
\end{itemize}
Again in this case, for any polynomial $P$ of critical degree, we have that $\Res_F\,(P)=  \pm
\frac{\det(\tilde\mm_P)}{\det(\tilde\mm)}.$
\end{example}
\medskip

\section{Toric Subresultants in the critical degree}\label{3}

In this section we define the toric subresultant of a monomial $h \in
S_{\rho}$ and show that this is precisely the numerator of the residue
of $h$. To set things up we must construct two complexes of free $\aa$
modules wich we will call the resultant and subresultant complexes
respectively. Along the way we will prove Theorem \ref{codimthm}
from the previous section. The approach will be to use Weyman's
complex \cite[Section 3.4E]{GKZ} to pass from an exact sequence of
sheaves to a generically exact complex of free modules.

\begin{proof}[Proof of Theorem \ref{codimthm}]
Note that, as $\Res_F\,(\Del) = \pm 1,$ $\Del$ is not in the ideal $\langle F_0, \dots, F_n \rangle$. Hence  the surjectivity of $\phi$ and the fact that $(S_F)_{\rho}$ has
codimension $1$ are equivalent.
\par
Let $F$ be our standard generic ample system. The polynomials $F_i$
are sections of sheaves $\mathcal{L}_i := \mathcal{O}(\alpha_i)$ on
$X$ (with coefficients in $\aa$). Given any subset $I$ of $\{0, \dots,
n\}$, let $\alpha_I = \sum_{i \in I} \alpha_i$. We get a corresponding
(dual) Koszul complex of sheaves:
\begin{equation}
0 \to \mathcal{O}(-\sum \alpha_i) \to \cdots \to
\underset{i}{\bigoplus} \mathcal{O}(-\alpha_i) \overset{F}{\to}
\mathcal{O}_X \to 0.
\end{equation}
If we now tensor this complex with the sheaf $\mathcal{M} :=
\mathcal{O}(\rho) = \mathcal{O}(\sum \alpha_i - \beta_0)$ we get a new complex:
\begin{equation} \label{sheafcom}
0 \to \mathcal{O}(-\beta_0) \to \cdots \to \underset{i}{\bigoplus} \mathcal{O}(\rho -\alpha_i) \overset{F}{\to} \mathcal{O}(\rho) \to 0.
\end{equation}
Since all sheaf Tor groups vanish when one of the factors is locally
free, it follows that the complex (\ref{sheafcom}) remains exact even
if $\mathcal{M}$ is not locally free.  We can therefore apply
``Weyman's complex'' \cite[Chapter 3, Theorem 4.11]{GKZ} which yields a double complex:
\begin{eqnarray*}
C^{-p, q} &=& H^q \left( X, \underset{0 \leq i_1 \cdots \leq i_p \leq n}{\bigoplus} \mathcal{L}_{i_1}^{\ast} \otimes \cdots \otimes \mathcal{L}_{i_p}^{\ast} \otimes \mathcal{M} \right)\\
          &=& \underset{|I| = n+1 - p}{\bigoplus} H^q(X, \OO(\alpha_I - \beta_0)).
\end{eqnarray*}
The corresponding total complex is generically exact with
differentials depending polynomially on the coefficients of the $F_i$;
therefore we can view this as a generically exact complex of free
$\aa$-modules. By the toric version of Kodaira vanishing, see \cite{Mus},
all cohomology terms in the complex vanish except when $q = 0$ or when
$p = n+1$ and $q = n$. Note that $H^0(X, \OO(\alpha_I - \beta_0)) =
S_{\alpha_I - \beta_0}$ and the differentials between these terms are
just those from the Koszul complex on $S$ determined by $F_0, \dots,
F_n$. Also the only non-vanishing higher cohomology term is $H^n(X,
\OO(-\beta_0)) \cong \CC$, which will correspond to a rank 1 free
$\aa$-module. The last differential $\phi \ : \ C^{-1} \to C^0$ can
therefore be chosen to be the map (\ref{modc}) from Section 3.
$$
\begin{array}{cccc}
\phi:&S_{\rho-\alpha_0}\oplus \dots \oplus S_{\rho-\alpha_n}\oplus\aa
&\to&S_\rho\\ &(G_0,\dots,G_n,c)&\mapsto&\sum_{i=0}^n G_iF_i+c\,\Del.
\end{array}
$$
This map is generically surjective and therefore we have proven \ref{codimthm}.
\end{proof}

\begin{definition}
The {\bf resultant complex} is the complex of free $\aa$-modules constructed
above. Namely
\begin{equation} \label{rescomp}
 0 \to S_{-\beta_0} \to \cdots \to \underset{i}{\bigoplus} S_{\rho -
 \alpha_i} \oplus \aa \overset{\phi}{\to} S_{\rho} \to 0
\end{equation}
where the map $\phi$ is as above.
\end{definition}

Let $\AA_{i}:=\{a\in\NN^s:\,\deg(a)=\alpha_i \}$, and $\ell$ the index of the
lattice spanned by $\cup \AA_{i}$ in $\mathbb{Z}^n$.

\begin{proposition} \label{resprop}
The determinant of the resultant complex (\ref{rescomp}) with respect
to the monomial bases is $$c \cdot {\rm res}_{\alpha_0, \dots,
\alpha_n}(F_0, \dots, F_n)^{\ell}$$ for some constant $c \in
\mathbb{Q}$.
\end{proposition}

\begin{proof} In the case when $\ell = 1$, so the $\AA_{i}$ span $\mathbb{Z}^n$, this is a consequence of Theorem 4.11 in \cite[Chapter 3.4E]{GKZ}. For the general case we note that the determinant still vanishes if and only if the resultant is 0, and the degree with respect to the coefficients of any $F_i$ is still
the mixed volume of all of the other supports with respect to the
given lattice $\mathbb{Z}^n$. The degree of the resultant, on the
other hand, is $\frac{1}{\ell}$ of this mixed volume \cite{PS}.
\end{proof}

\begin{conjecture}
The constant $c = \pm 1$.
\end{conjecture}

The following corollary also appears, under weaker hypothesis on $F$, in
\cite{CoD}. 

\begin{corollary}
The complex below of free $\aa$-modules is generically exact
everywhere but at the last step.

\begin{equation} \label{comp}
 0 \to S_{-\beta_0} \to \cdots \to \underset{i}{\bigoplus} S_{\rho -
 \alpha_i} \overset{F}{\to} S_{\rho}.
\end{equation}

\end{corollary}

There are two ways to enforce exactness at the last stage of the
complex.  The resultant complex does this by enlarging the second to
last module.  The second way to get an exact complex out of
(\ref{comp}) is to corestrict the last map to a smaller target. Pick a
monomial $h \in S_{\rho}$ and define $S_{\rho}/h := \aa \langle x^a,
\deg(a) = \rho, x^a \neq h \rangle$.

\begin{definition} The {\bf subresultant complex} with respect to $h$ is the complex
\begin{equation} \label{subrescomp}
 0 \to S_{-\beta_0} \to \cdots \to \underset{i}{\bigoplus} S_{\rho -
 \alpha_i} \overset{F_h}{\to} S_{\rho}/h \to 0
\end{equation}
Here $F_h$ is the multiplication map of the $F_i$ corestricted to $S_{\rho}/h$.
\end{definition}

\begin{proposition} If $h$ does not belong to the ideal $\langle F_0, \dots, F_n\rangle$, with coefficients in $Q(\aa)$, then the complex (\ref{subrescomp}) is generically exact.
\end{proposition}

\begin{proof} This is an immediate consequence of Theorem \ref{codimthm}. \end{proof}

\par Now, as the homology with coefficients in $\aa$ vanishes for
$p>0,$ the determinant of the complex with respect to the monomial
bases is an element of $\aa$ (\cite[Appendix A]{GKZ}) and may be
computed as the $\gcd$ of the maximal minors of $F_h.$ So we can
define:

\begin{definition} \label{subresdef}
If $h$ generates ${S_F}_\rho,$ then let
$$\ss_h:=\det(\mbox{complex}\,(\ref{subrescomp}))\in \aa,$$ where the
determinant is taken with respect to the monomial bases of $\KK.$ If
$h$ belongs to ideal generated by $F_0,\dots,F_n,$ then we set
$\ss_h:=0.$ The polynomial $\ss_h$ is well-defined, up to a sign, and
is called the $h$-subresultant of the family (\ref{generic}).
\end{definition}

\smallskip
\begin{proposition} $_{}$
\par
\begin{enumerate}
\item For each $i$ and each monomials $h$ of degree $\rho,$ $\ss_h$
is homogeneous in the coefficients of $F_i.$ If it is not identically
zero, it has total degree equal to $\sum_i MV(\alpha_0, \dots, \alpha_{i-1},
\alpha_{i+1}, \dots, \alpha_n) -1.$ Here $MV$ is the mixed volume of
the polytopes corresponding to the given ample degrees.

\item Let ${\bf k}$ be a field of characteristic zero. For every
specialization of the coefficients of $F_i$ in ${\bf k},$ we have that
$\ss_h\neq0$ if and only if
$${\bf k}\langle h\rangle +\langle F_0, \dots F_n \rangle_\rho={\bf k}[x_0,\dots,x_{s-1}]_\rho.$$
\end{enumerate}
\end{proposition}

\smallskip
\begin{proof}
The second statement is a straightforward consequence of the
definition of $\ss_h$ as the determinant of the complex
(\ref{subrescomp}).  The first part follows by comparing the determinant
of (\ref{subrescomp}) and the determinant of the resultant complex
(\ref{rescomp}), whose degree in the coefficients of $f_i$ equals
$MV(\alpha_0, \dots, \alpha_{i-1}, \alpha_{i+1}, \dots, \alpha_n)$.
\end{proof}

\section{Residues, Resultants and Subresultants}\label{4}

We are now ready for our second main theorem.

\begin{theorem}\label{mtt}
$$\Res_F\,(h)=\pm{\frac{\ss_h}{c\,.\,{\rm res}_{\alpha_0,\dots,
\alpha_n}(F_0,\dots,F_n)^\ell}}.$$
\end{theorem}

\begin{proof}

Using the result of Theorem \ref{quot}:
$$\Res_F\,(h)=\pm{\frac{\det(\tilde{\mm_h})}{\det(\tilde{\mm})}} =
\pm{\frac{\delta_1\ss_h}{\delta_1 \,c\,{\rm
res}_{\alpha_0,\dots,\alpha_n}(F_0,\dots,F_n)^\ell}}.$$
This is due to
the fact that the extraneous factor $\delta_1$ for a maximal minor in
the resultant complex and for the corresponding maximal minor in any
subresultant complex are the same. The result now follows from the
definition of the subresultant and Proposition \ref{resprop}.
\end{proof}

\begin{corollary}
We get the following factorization in $\aa:$
$$\ss_h={\rm
res}_{\alpha_0,\dots,\alpha_n}(F_0,\dots,F_n)^{\ell-1}\,P_h,$$ where $P_h$
is a polynomial which is not a factor of the resultant.
\end{corollary}
\begin{proof}
In \cite[Theorem 1.4]{CDS} it is shown that the residue is a rational
function whose denominator is ${\rm
res}_{\alpha_0,\dots,\alpha_n}(F_0,\dots,F_n).$ While their proof is in the
generalized unmixed setting, it carries over to our setting as well,
since it relies only on the representation of the toric residue as a
sum of local residues, \cite[Theorem 0.4]{CCD}.  Counting degrees, it
turns out that $P_h$ has degree in the coefficients of $f_i$ one less
than the sparse resultant. So, it cannot be a factor of it.
\end{proof}

\medskip
Theorem \ref{mtt} may be regarded as a generalization of Jouanolou's results in the dense case. Suppose that $F_0,\dots,F_n$ are generic homogeneous
polynomials of respective degrees $d_0,\dots,d_n.$ In \cite[(2.9.6)]{J1}, a linear function 
$\omega:\left(S / \langle F_0,\dots,F_n\rangle\right)_\rho\to \aa$ is defined by setting $\omega(\Del):=res_{d_0,\dots,d_n}(F_0,\dots,F_n).$ Hence, $\omega(h)$ may be
regarded as the numerator of the residue of $h.$ Several properties of this morphism are studied in \cite{J2,J3}. In \cite[Corollaire 3.9.7.7]{J3}, it is
shown that $\omega(h)$ may be computed as a quotient of two determinants. Comparing this quotient with Chardin's recipe for computing the subresultant as a
quotient of two determinants (\cite{Cha2}) we get that, if $h$ is a monomial then $\omega(h)$ is the classical subresultant of the
set $\{h\}$ with respect to $F_0,\dots,F_n.$
\medskip

\begin{example}
We present here an example where $\ell>1.$ Let $P$ be the simplex in $\RR^3$ which is the convex hull of $(0,0,0),\,(0,1,0),\,(0,0,1)$
and $(3,1,1).$ In this case, $\ell=3$ and $S$ is a ring of polynomials in $4$ variables. Let $\alpha=P\cap\ZZ^3,$ and consider the following four generic polynomials in $S_\alpha:$
$$F_i:=a_ix_1^3+b_ix_2^3+c_ix_3^3+d_ix_4^3,\ i=0,1,2,3.$$
$\Delta_{\sigma}$  of this system equals $x_1^2x_2^2x_3^2x_4^2$ times the determinant of
$$D:=\left(\begin{array}{cccc}
a_0&b_0&c_0&d_0\\
a_1&b_1&c_1&d_1\\
a_2&b_2&c_2&d_2\\
a_3&b_3&c_3&d_3
\end{array}\right).$$
Also, it is easy to see that ${\rm res}_{\alpha,\alpha,\alpha,\alpha}(F_0,F_1,F_2,F_3)=\det(D).$
The matrix of the last morphism of the complex (\ref{modc}) has size $33\times21$ in this case. A nonzero maximal minor of this matrix is
{\tiny $$\left(\begin{array}{ccccccccccccccccccccc}
0&a_1&b_1&c_1&d_1&0&0&0&0&0&0&0&0&0&0&0&0&0&0&0&0\\
0&a_2&b_2&c_2&d_2&0&0&0&0&0&0&0&0&0&0&0&0&0&0&0&0\\
0&a_3&b_3&c_3&d_3&0&0&0&0&0&0&0&0&0&0&0&0&0&0&0&0\\
0&a_4&b_4&c_4&d_4&0&0&0&0&0&0&0&0&0&0&0&0&0&0&0&0\\
0&0&0&0&a_1&b_1&c_1&d_1&0&0&0&0&0&0&0&0&0&0&0&0&0\\
0&0&0&0&a_2&b_2&c_2&d_2&0&0&0&0&0&0&0&0&0&0&0&0&0\\
0&0&0&0&a_3&b_3&c_3&d_3&0&0&0&0&0&0&0&0&0&0&0&0&0\\
0&0&0&a_1&0&0&d_1&0&b_1&0&0&0&0&0&0&0&0&0&0&0&c_1\\
0&0&0&a_2&0&0&d_2&0&b_2&0&0&0&0&0&0&0&0&0&0&0&c_2\\
0&0&0&0&0&0&0&0&0&a_1&b_1&c_1&d_1&0&0&0&0&0&0&0&0\\
0&0&0&0&0&0&0&0&0&a_2&b_2&c_2&d_2&0&0&0&0&0&0&0&0\\
0&0&0&0&0&0&0&0&0&a_3&b_3&c_3&d_3&0&0&0&0&0&0&0&0\\
0&0&0&0&0&0&0&0&0&a_4&b_4&c_4&d_4&0&0&0&0&0&0&0&0\\
0&0&0&0&0&0&0&0&0&0&0&0&a_1& b_1&d_1&0&0&c_1&0&0&0\\
0&0&0&0&0&0&0&0&0&0&0&0&a_2&b_2&d_2&0&0&c_2&0&0&0\\
0&0&0&0&0&0&0&0&0&0&0&0&a_3&b_3&d_3&0&0&c_3&0&0&0\\
0&0&0&0&0&0&0&0&0&0&0&a_1&0&0&0&b_1&c_1&d_1&0&0&0\\
0&0&0&0&0&0&0&0&0&0&0&a_2&0&0&0&b_2&c_2&d_2&0&0&0\\
0&0&a_1&0&0&d_1&0&0&c_1&0&0&0&0&0&0&0&0&0&b_1&0&0\\
0&0&0&0&0&0&0&0&0&0&a_1&0&0&d_1&0&c_1&0&0&0&b_1& 0\\
\det(D)&0&0&0&0&0&0&0&0&0&0&0&0&0&0&0&0&0&0&0& 0
\end{array}
\right).$$ }
Its determinant equals
$${\det(D)}^3\,
(c_2b_1- b_2c_1)^2(d_3c_2b_1 - b_2c_1d_3 - d_1b_3c_2 - c_3d_2b_1 + c_3d_1b_2 + d_2b_3c_1)^2
$$
so we have that
$$\delta_1=(c_2b_1- b_2c_1)^2(d_3c_2b_1 - b_2c_1d_3 - d_1b_3c_2 - c_3d_2b_1 + c_3d_1b_2 + d_2b_3c_1)^2.$$
The first column of this matrix is indexed by $x_1^2x_2^2x_3^2x_4^2.$ It is easy to see that, for every monomial $h$ of critical degree, $\ss_h=0$ unless $h=x_1^2x_2^2x_3^2x_4^2.$
In this case, $\ss_{x_1^2x_2^2x_3^2x_4^2}={\det(D)}^2$ and hence $P_{x_1^2x_2^2x_3^2x_4^2}=\pm1.$
\par This can be explained as follows: every monomial of critical degree is a multiple of $x_i^3$ for at least one $i=0,\dots,3$ except
$x_1^2x_2^2x_3^2x_4^2.$ An easy consequence of Cramer's rule is that every monomial which is multiple of $x_i^3$ is in the ideal generated by the
generic $F_i$'s. This is why all except one of the subresultants are identically zero.
\end{example}

\medskip
\begin{conjecture}
If $P_h$ is not identically zero, then it is an irreducible element of $\aa.$ In particular, when $\ell=1,$ every subresultant $\ss_h$ is irreducible.
\end{conjecture}
\medskip

\section{Computing global residues ``\'a la Macaulay''}\label{smac}
In this section, we will review the toric algorithm of \cite{CD} for computing global residues by means of toric residues (see also \cite{CDS}). As a straightforward consequence of
their algorithm and our results,
we get a quotient formula for computing global residues. In the dense case, we recover the quotient type formula given by Macaulay in \cite{Mac} for computing the global residue of
$x_1^{d_1-1}x_2^{d_2-1}\dots x_n^{d_n-1}$ with respect to a generic family of polynomials of degrees $d_1,\dots,d_n.$
\par
Let ${\bf A}_1,\dots,{\bf A}_n$ subsets of $\ZZ^n,$ and consider $n$ Laurent polynomials in $n$ variables $t_1,\dots,t_n$ having support in ${\bf A}_1,\dots,{\bf A}_n$ respectively:
$$ f_j \quad = \, \sum_{m\in {\bf A}_j}
u_{jm}\,\cdot\,t^m\ j=1,\dots,n.$$
Let $V$ be the  set of common zeros of
$f_1,\ldots,f_n$ in the
torus $T = \left( \CC^* \right)^n $. If $V$ is finite and all its roots are simple, then
for any Laurent polynomial
$\,q\in\CC[t_1^{\pm 1}, \ldots,t_n^{\pm 1}]$,  the
{\it global residue} of the differential form
$$ \phi_q \quad = \quad \frac q{f_1\cdots f_n}\,
\frac {dt_1}{t_1}\wedge
\cdots\wedge \frac {dt_n}{t_n}\,, $$
is defined as
$\sum_{\xi \in V} \frac {q(\xi)}{J^T(f)(\xi)},$
where $J^T(f)$ denotes the {\it affine toric Jacobian\/}
$$J^T(f) \quad :=  \quad
\det \bigl( t_k\,\frac {\partial f_j} {\partial t_k}\bigr)_{1\leq
j,k \leq n}\, .$$
Global residues are basic invariants of
multivariate polynomial systems (see \cite{CD,CDS} and the references therein).
\par The link between toric and global residues is given in \cite[Theorem 4]{CD}:
\begin{theorem}
Let $f_1,\dots,f_n\in\CC[t_1^{\pm1}\,\dots,t_n^{\pm1}]$ be Laurent
polynomials having a finite number of zeroes in $T$, and $g$ another
Laurent monomial. Then there is a projective toric variety $X$ with
homogeneous coordinate ring $S_X$ and a homogeneous element $F_0\in
S_X$ such that
\begin{enumerate}
\item There is an ``homogenization rule'' which assigns to every $f_i$ a homogeneous polynomial $F_i\in S_X,\,i=1,\dots,n;$
\item the family $F_0,\dots,F_n$ has no zeros in $X;$
\item There is another homogeneous monomial $G\in S_X$ such that
$$\mbox{\bf Global Residue}\,_f\,(g)= \Res_F\,(G).$$
\end{enumerate}
\end{theorem}
\medskip
As an immediate consequence, we get also a quotient formula for computing global residues, as the following example shows:

\begin{example} \label{ddo}
This example has already appeared in the introduction of \cite{CDS}.
We want to compute the global residue of $g(t):=t_1^3t_2^2$ with respect to the generic system
$$ \begin{array}{ccl}
f_1 & =  &
 a_0  t_1^2 \,+\, a_1  t_1 t_2 \,+\, a_2  t_2^2 \,+\, a_3  t_1
\,+\, a_4  t_2 \,+\, a_5\,,
 \\
f_2 & = &
b_0 t_1^2 \,+\, b_1  t_1 t_2 \,+\, b_2  t_2^2 \,+\, b_3  t_1
\,+\, b_4  t_2 \,+\, b_5\,.
\end{array}
$$
Applying the algorithm of \cite{CD} we get that, in this case,
the corresponding toric variety $X$ is $\mathbb{P}^2$ with the
standard homogeneous coordinates $S_X=\mathbb{C}[x_0, x_1, x_2].$ The ``homogeneous'' polynomials are
$$\begin{array}{l}
F_0(x_0,x_1,x_2):=x_0^2\\
F_1(x_0,x_1,x_2):=a_0x_1^2+a_1x_1x_2+a_2x_2^2+a_3x_1x_0+a_4x_2x_0+a_5x_0^2\\
F_2(x_0,x_1,x_2):=b_0x_1^2+b_1x_1x_2+b_2x_2^2+b_3x_1x_0+b_4x_2x_0+b_5x_0^2,
\end{array}$$
and $G=x_1^2x_2.$
\par The global residue may be computed then as the toric residue of $G$ with respect to $F_0,F_1,F_2.$ Using our methods, it turns out that 
the matrix $\MM$ of Theorem \ref{quot} is square. Computing it explicitly, we get
$$
\MM=\bordermatrix{
&x_0x_1^2&x_0x_1x_2&x_0x_2^2&x_0^2x_1&x_0^2x_2&x_0^3&x_1^3&x_1^2x_2&x_1x_2^2&x_2^3\cr
&a_0&a_1&a_2&a_3&a_4&a_5&0&{\bf 0}&0&0\cr
&b_0&b_1&b_2&b_3&b_4&b_5&0&{\bf 0}&0&0\cr
&0&0&0&0&0&1&0&{\bf 0}&0&0\cr
&a_3&a_4&0&a_5&0&0&a_0&{\bf a_1}&a_2&0\cr
&b_3&b_4&0&b_5&0&0&b_0&{\bf b_1}&b_2&0\cr
&0&0&0&1&0&0&0&{\bf 0}&0&0\cr
&0&a_3&a_4&0&a_5&0&0&{\bf a_0}&a_1&a_2\cr
&0&b_3&b_4&0&b_5&0&0&{\bf b_0}&b_1&b_2\cr
&0&0&0&0&1&0&0&{\bf 0}&0&0\cr
&{\bf j_1}&{\bf j_2}&{\bf j_3}&{\bf j_4}&{\bf j_5}&{\bf j_6}&{\bf j_7}&{\bf j_8}&{\bf j_9}&{\bf j_{10}}
},
$$
where the ${\bf j_k}$ are the coefficients of $\Del$.
The matrix $\tilde\MM_h$ is made by deleting in $\mm$ the eighth row and the last column.
Hence, we have in virtue of Theorem \ref{quot}
$$\Res_F\,(h)=\pm \frac{\det(\tilde\MM_h)}{\det(\MM)}=\pm\frac{\det(\tilde\MM_h)}{{\rm res}_{2,2,2}(F_0,F_1,F_2)},$$
and one can check that $\det(\tilde\MM_h)$ is the polynomial $P_{32}$ of the introduction of \cite{CDS}.
\end{example}

\bigskip
We close this section by showing that the method for computing global residues as a quotient of two determinants presented here, may be regarded as a generalization of a formula given by Macaulay in the classical case.
In order to follow his notation, let
$F_1(x_0,\dots,x_n),\dots, \\ F_n(x_0,\dots,x_n)$ be generic homogeneous polynomials of respective degrees $d_1,\dots,d_n,$ and
set $f_i:=F_i(1,x_1,\dots,x_n).$ Let $J$ be the affine Jacobian of the $f_i,$ i.e. $J=\det\left(\frac{\partial f_i}{\partial x_j}\right)_{i,j},$ and
$$V(f_1,\dots,f_n):=\{\xi_1,\dots,\xi_{d_1\dots d_n}\} \\ \subset\overline{Q(\aa)}$$ be the variety defined by the common zeroes of the $f_i$ in the algebraic
closure of $Q(\aa).$ We denote with ${\bf m}$ the monomial $x_1^{d_1-1}\dots x_n^{d_n-1}.$
\par From display $(13)$ in \cite{Mac} and following his notation, we have
\begin{equation}\label{gmacaulay}
\sum_{j=1}^{d_1\dots d_n}\frac{{\bf m}(\xi_j)}{J(\xi_j)}=\pm \frac{R(n,t_{n}-1)}{R(n,t_{n})},
\end{equation}
where
\begin{itemize}
\item $R(n,t_n)$ is the resultant of $F_1(0,x_1,\dots,x_n),\dots,F_n(0,x_1,\dots,x_n),$
\item $R(n,t_n-1)$ is the subresultant of the monomial ${\bf m}$ with respect to
$F_1(0,x_1,\dots,x_n),\dots,F_n(0,x_1,\dots,x_n).$
\end{itemize}
Comparing the left hand side of (\ref{gmacaulay}) with the definition of the definition given above, we have that 
(\ref{gmacaulay}) is actually the global residue of $x_1^{d_1}\dots x_n^{d_n}$ with respect to $f_1,\dots,f_n$.
Applying the toric algorithm of \cite{CD} we get the following:
\begin{itemize}
\item $X=\PP^n,\,S_X=\CC[x_0,\dots,x_n]$ with homogenization given by total degree;
\item $F_0=x_0,$ and if $g=x_1^{d_1}\dots x_n^{d_n},$ then $G={\bf m}.$
\end{itemize}
Denote with ${\rm res}$ the homogeneous resultant for a family of $n+1$ homogeneous polynomials in $n+1$ variables
of degrees $1,d_1,\dots d_n.$
Then, it turns out that (\ref{gmacaulay}) equals $\Res_{F}\,({\bf m})$.
Applying Theorem \ref{mtt}, we can write $\Res_F\,({\bf m})$ as $\pm\frac{\ss_{\bf m}}{{\rm res}(x_0,F_1,\dots,F_n)}.$
Now, specializing a generic $F_0$ to $x_0$ and
applying \cite[Lemma 1]{Cha2}, we get that $\ss_{\bf m}\mapsto R(n,t_n-1)$ and ${\rm res}\mapsto R(n,t_n).$

\medskip
\subsection{Acknowledgements} We would like to thank Bernd Sturmfels for his insightful comments on preliminary versions of this paper. We thank Eduardo Cattani and Alicia Dickenstein for pointing out that the original proof of Theorem 2.1 was incomplete and providing a suggestion on how to fill the gap. We are also grateful to David Cox for helpful discussions and along with Alicia Dickenstein for having communicated to us their work in progress.  All the computations
were made with {\tt Maple}.  The first author was supported by the
Miller Institute for Basic Research in Science, in the form of a
Miller Research Fellowship (2002--2005).


\begin{thebibliography}{XXX}

\bibitem[Bez]{Bez}
\'E.~B\'ezout.
\newblock{\em Th\'eorie g\'en\'erale des \'equations alg\'ebriques.}
\newblock Ph.-D. Pierres: Paris.

\bibitem[BC]{BC}
V.~Batyrev and D.~Cox.
\newblock{\em On the Hodge structure of projective hypersurfaces in toric varieties,}
\newblock Duke J. Math. 75 (1994), 293-338

\bibitem[BM]{BM}
V.~Batyrev and D.~Cox.
\newblock{\em Mixed toric residues and Calabi-Yau complete intersections.}
\newblock math.AG/0206057

\bibitem[Cha1]{Cha1}
M.~Chardin.
\newblock{\em Multivariate subresultants.}
\newblock J. Pure Appl. Algebra 101, No.2, 129-138 (1995).

\bibitem[Cha2]{Cha2}
M.~Chardin.
\newblock{\em Formules \`a la Macaulay pour les sous-r\'esultants en
plusieurs variables.}
\newblock C. R. Acad. Sci. Paris S\'er. I Math.  319  (1994),  no. 5, 433--436.

\bibitem[CCD]{CCD}
E.~Cattani, D.~Cox, and A.~Dickenstein.
\newblock{\em Residues in toric varieties.}
\newblock Compositio Math. 108 (1997), no. 1, 35--76.

\bibitem[CoD]{CoD}
D.~Cox and A.~Dickenstein.
\newblock{\em Vanishing and codimension one theorems on complete toric varieties}
\newblock In progress, personal communication

\bibitem[CD]{CD}
E.~Cattani and A.~Dickenstein.
\newblock{\em A global view of residues in the torus.}
\newblock Algorithms for algebra (Eindhoven, 1996).
J. Pure Appl. Algebra 117/118 (1997), 119--144.

\bibitem[CDS1]{CDS1}
E.~Cattani, A.~Dickenstein, and B.~Sturmfels.
\newblock{\em Computing multidimensional residues.}
\newblock Algorithms in algebraic geometry and applications (Santander, 1994),  135--164, Progr. Math., 143, Birkh\"auser, Basel, 1996.

\bibitem[CDS2]{CDS}
E.~Cattani, A.~Dickenstein, and B.~Sturmfels.
\newblock{\em Residues and resultants.}
\newblock J. Math. Sci. Univ. Tokyo 5 (1998), no. 1, 119--148.

\bibitem[Cox1]{cox1}
D.~Cox.
\newblock{\em The homogeneous coordinate ring of a toric variety.}
\newblock J. Algebr. Geom., 4, 17--50 (1995).

\bibitem[Cox2]{cox}
D.~Cox.
\newblock{\em Toric residues.}
\newblock Ark. Mat. 34 (1996), no. 1, 73--96.

\bibitem[CLO1]{CLO1}
D.~Cox, J.~Little, and D.~O'Shea.
\newblock{\em Ideals, Varieties, and Algorithms.
An introduction to computational algebraic geometry and commutative algebra.}
\newblock Undergraduate Texts in Mathematics. New York, NY: Springer. xiii, 1996.

\bibitem[CLO2]{CLO2}
D.~Cox, J.~Little, and D.~O'Shea.
\newblock{\em Using Algebraic Geometry}.
\newblock  Graduate Texts in Mathematics, 185. Springer-Verlag, New York, 1998.

\bibitem[Ful]{ful}
W.~Fulton
\newblock{\em Introduction to Toric Varieties.\/}
\newblock Ann. of Math. Studies, no. 131, Princeton NJ, Princeton University Press (1993).

\bibitem[GH]{GH}
P.~Griffiths and J.~Harris.
\newblock{\em Principles of Algebraic Geometry}.
\newblock John Wiley \& Sons, New York, 1978.

\bibitem[GK]{GK} 
O.A.~Gelfond, and A.G.~Khovanskii. 
\newblock{\em Toric geometry and Grothendieck residues.\/}  
\newblock Mosc. Math. J.  2  (2002),  no. 1, 99--112, 199. 14M25

\bibitem[GKZ]{GKZ}
I.M.~ Gelfand, M.M.~ Kapranov, and A.V.~ Zelevinsky.
\newblock{\em Discriminants, Resultants, and Multidimensional Determinants.}
\newblock Mathematics: Theory \& Applications. Birkh\"auser Boston, Inc., Boston, MA, 1994.

\bibitem[J1]{J1}
J.P.~Jouanolou. 
\newblock{\em Singularit\'es rationnelles du r\'esultant.\/}  
\newblock Algebraic geometry (Proc. Summer Meeting, Univ. Copenhagen, Copenhagen, 1978),  
pp. 183--213, Lecture Notes in Math., 732, Springer, Berlin, 1979. 

\bibitem[J2]{J2}
J.P.~Jouanolou. 
\newblock {\em Aspects invariants de l'élimination.} 
\newblock Adv. Math.  114  (1995),  no. 1, 1--174. 

\bibitem[J3]{J3}
J.P.~Jouanolou.
\newblock{ \em Formes d'inertie et r\'esultant: un formulaire.} (French) 
\newblock Adv. Math.  126  (1997),  no. 2, 119--250. 

\bibitem[Kho]{Kho}
A.G.~Khovanskii.
\newblock{\em Newton polyhedra and toroidal varieties.\/}
\newblock Russian Math. Surveys 33 (1978), 237--238.

\bibitem[Mac]{Mac}
F.~Macaulay.
\newblock{\em Some formulae in elimination}.
\newblock Proc. London Math. Soc. 1 \textbf{33}, 3--27, 1902.

\bibitem[Mus]{Mus}
M.~Musta{\c{t}}{\u{a}}
\newblock{\em Vanishing theorems on toric varieties}.
\newblock Tohoku Math. J.(2) \textbf{54} (2002), no. 3, 451--470.

\bibitem[PS]{PS}
P.~Pedersen and B.~Sturmfels,
\newblock {\em Product formulas for resultants and Chow forms}.
\newblock Mathematische Zeitschrift \textbf{214} (1993), 377-396.

\bibitem[Sop]{Sop}
I.~Soprounov.
\newblock{\em Residues and tame symbols in toric geometry.\/}
\newblock math.AG/0203114.

\bibitem[Sza]{sza}
 A.~ Szanto.
\newblock{\em  Multivariate subresultants using Jouanolou's resultant matrices.}
\newblock Preprint.

\bibitem[VdW]{VdW}
B.~L.~van der Waerden.
\newblock{\em Modern Algebra.\/}
\newblock $3^{rd}$ edn. New York, F.~Ungar Publishing Co., 1950.
\end{thebibliography}
\end{document}